\documentclass[12pt,usenames,dvipsnames]{article}
\usepackage[a4paper]{anysize}\marginsize{2cm}{2cm}{2cm}{2cm}
\pdfpagewidth=\paperwidth \pdfpageheight=\paperheight
\usepackage{amsfonts,amssymb,amsthm,amsmath,eucal,tabu,hyperref}
\usepackage{mathabx}
\usepackage{pgf}
\usepackage{array}
\usepackage{pstricks}
\usepackage{tikz-cd}
\usepackage{pstricks-add}
\usepackage{booktabs} 
\usepackage{multicol} 
\usepackage{pgf,tikz,pgfplots}
\usepackage{float}
\pgfplotsset{compat=1.10}
\usetikzlibrary{automata}
\usetikzlibrary{arrows}
\usetikzlibrary{matrix}
\usepackage{indentfirst}
\usepackage{xargs} 
\usepackage{xcolor}
\usepackage{todonotes}

\theoremstyle{plain}
\newtheorem{thm}{Theorem}[section]
\newtheorem{theorem}[thm]{Theorem}
\newtheorem{lemma}[thm]{Lemma}
\newtheorem{proposition}[thm]{Proposition}

\newtheorem{corollary}[thm]{Corollary}

\theoremstyle{definition}
\newtheorem{definition}[thm]{Definition}
\newtheorem{remark}[thm]{Remark}
\newtheorem{example}[thm]{Example}

\newtheorem{thevarthm}[thm]{\varthmname}

\newenvironment{varthm*}[1]{\trivlist\item[]{\bf #1.}\it}{\endtrivlist}

\newcommand{\dP}{\check{\mathbb{P}}}

\def\keywordname{{\bfseries Keywords}}%
\def\keywords#1{\par\addvspace\medskipamount{\rightskip=0pt plus1cm
		\def\and{\ifhmode\unskip\nobreak\fi\ $\cdot$
		}\noindent\keywordname\enspace\ignorespaces#1\par}}
\def\subclassname{{\bfseries Mathematics Subject Classification
		(2010)}\enspace}
\def\subclass#1{\par\addvspace\medskipamount{\rightskip=0pt plus1cm
		\def\and{\ifhmode\unskip\nobreak\fi\ $\cdot$
		}\noindent\subclassname\ignorespaces#1\par}}

\def\P{\mathbb{P}}
\def\A{\mathcal{A}}
\def\K{\mathbb{K}}
\def\C{\mathbb{C}}

\def\N{\mathbb{N}}

\def\O{\mathcal{O}}
\def\eL{\mathcal{L}}
\def\H{\mathcal{H}}
\DeclareMathOperator{\Der}{Der}

\DeclareMathOperator{\D}{D}
\DeclareMathOperator{\syz}{syz}
\DeclareMathOperator{\sym}{Sym}

\DeclareMathOperator{\adim}{adim}
\DeclareMathOperator{\vdim}{vdim}

\title{Unexpected hypersurfaces of type $(d+k,d)$
}
\author{Marek Janasz, Grzegorz Malara, Halszka Tutaj-Gasi\'nska}
\date{\today}
\begin{document}
	\maketitle
	\begin{abstract} Unexpected hypersurfaces arise when  vanishing in points of a set $Z$ and higher-order vanishing along a general linear subspace fails to impose the expected number of independent conditions on forms of a fixed degree. The phenomenon was first observed for planar curves by Cook, Harbourne, Migliore and Nagel.
    This paper shows a syzygy-based construction of, possibly unexpected, hypersurfaces of degree $d+k$ in $\P^n$, vanishing along a codimension two general linear subspace with multiplicity $d$;
    thus generalizing the work of Trok and the previous work of the last two authors. 
    Our framework unifies the classical planar cases with higher-dimensional examples, including Trok’s construction. We give a sufficient  criterion for unexpectedness (via the splitting behaviour the syzygy bundles of the powers of the Jacobian ideal,
    associated with the hyperplane arrangement dual to 
    $Z$) and provide explicit examples in  $\P^3$ and $\P^4$.

		\subclass{14N20, 14C20, 13D02, 05E40, 14F05}

        
	\end{abstract}
\maketitle
	
\section{Introduction}

The notion of an \emph{unexpected hypersurface} was introduced in \cite{CHMN}, where the authors, using results of \cite{FV,V}, presented a quartic curve in $\P^2$ passing through a special set $Z$ of nine points that impose independent conditions on quartics, and having a triple point in a general position $B$. A simple parameter count shows that such a curve should not exist, as the imposed conditions exactly exhaust the expected dimension of the space of quartics. Its existence therefore constitutes a genuine geometric anomaly, and the curve is called \emph{unexpected} of type $(4,3)$. Since that seminal example, the concept of unexpected curves has been generalised to \emph{unexpected hypersurfaces}; see, for instance, the introduction to \cite{HMT-G}. Formally, let $\K[x_0,\ldots,x_n]$ be the homogeneous coordinate ring of $\P^n=\P^n_{\K}$, where $\K$ is an algebraically closed field of characteristic $0$. For a subscheme $Z\subset\P^n$ with saturated ideal $I_Z$, and general linear subvarieties $B_i\subset\P^n$ of multiplicities $m_i>0$, we set $B=m_1B_1+\cdots+m_rB_r$ and denote by $H_B$ the Hilbert polynomial of $\K[x_0,\ldots,x_n]/I_B$.
\begin{definition}[{\cite{HMT-G}}]
A hypersurface $H$ of degree $d$ passing through $Z$ and vanishing along each $B_i$ with multiplicity $m_i$ is said to be \emph{unexpected} if
$$\dim [I_{Z\cup B}]_d>\max\left\{0,\dim[I_Z]_d-H_B(d)\right\},$$
that is, if vanishing on $B$ imposes fewer than the expected number of conditions on the forms of degree $d$.
\end{definition}
In the original setting, it is typically assumed that the conditions imposed by $Z$ are independent. In our approach, we relax this assumption: the set $Z$ may impose dependent conditions on forms of degree $d$, but the degree of dependence is explicitly taken into account when computing the virtual dimension $\vdim$. Whenever the actual dimension still satisfies $\adim>\vdim$, the corresponding system remains of interest and is regarded as exhibiting unexpected behaviour.

After \cite{CHMN}, numerous examples of unexpected varieties have been found, mostly in $\P^2$, and described from different perspectives, see, e.g., \cite{BMSS, DIV, DMO, Dimca, DFHMST,DHRST,DMTG25,FGST21,FGHM, HMNT,HMT-G, KS,MT-G,S, Trok}.  
In most cases only one general linear subvariety of dimension $0$ (a point) is considered, 
so that the unexpectedness arises when adding a fat point $mB$ makes the system dependent, 
although $Z$ itself imposes independent conditions.  
In the most interesting cases, the resulting hypersurface is irreducible and its existence is far from obvious.  
Known proofs of existence typically rely either on vector bundle methods (especially on the splitting type 
of logarithmic or syzygy bundles) or on direct computer computation of the defining equation. A different approach to proving the existence of unexpected hypersurfaces was proposed in \cite{DFHMST} and another one in \cite{DMTG25}, where the authors provided a combinatorial proof of this phenomenon. Beyond the rapidly growing catalogue of examples in the plane, a vector-bundle viewpoint has proved particularly fruitful; see, for instance, the analysis of logarithmic bundles and splitting types in \cite{FV, Trok}. At the same time, matrix and syzygy methods have provided constructive proofs and effective checks for existence in concrete configurations \cite{DMO,DFHMST,MT-G}. The present work merges these perspectives in higher dimension: we keep the syzygetic backbone, but replace the general point by a general codimension-two subspace $\H\subset\P^n$, and we pass from first-order relations to $k$-th order derivations. This yields a uniform mechanism that produces hypersurfaces of type $(d+k,d)$ and expresses their unexpectedness in terms of the splitting data of $D_0^k(\A_Z)$ restricted to the dual line $\eL=\check\H$.

The present paper continues the line of research initiated in~\cite{MT-G}, extending the syzygy-based construction of unexpected curves in~$\P^2$ to higher-dimensional projective spaces. In contrast to the planar case, the role of a general point is replaced by a general codimension-two subspace $\H\subset\P^n$. This shift allows to treat unexpected hypersurfaces of type $(d+k,d)$ in arbitrary dimension $n$, and to relate their existence to the behaviour of higher-order derivation bundles $D^k(\A_Z)$. Geometrically, our approach gives a unified interpretation of both the classical planar phenomena and Trok’s higher-dimensional construction \cite{Trok}, showing that the same syzygetic mechanism governs unexpectedness across dimensions.

The structure of the paper is as follows. Section \ref{sec:basic} recalls basic definitions and notation concerning dual projective spaces, hyperplane arrangements, and derivation modules. In Section \ref{sec:construction-generalized} we develop a general syzygy-based construction producing hypersurfaces of degree $d+k$ with prescribed multiplicity along a codimension-two subspace. Section \ref{sec:derivations} introduces the framework of $k$-derivations and establishes an explicit correspondence between the restricted bundles $D^k(\A_Z)\big|_{\eL}$ and the graded modules $I(Z,\H)(-k)$, generalizing the planar results of \cite{MT-G} and the constructions of \cite{Trok}. Section \ref{sec:unexpectedness} provides a sufficient numerical criterion for unexpectedness in terms of the splitting type of the bundle $D_0^k(\A_Z)$, together with a compact closed formula for the expected dimension of the relevant linear systems. Finally, in Section \ref{sec:examples} we present explicit examples in $\P^4$ and $\P^3$, computed in \texttt{Singular}, which illustrate the theory and demonstrate the occurrence of unexpected hypersurfaces arising from crystallographic configurations.

\section{Basic facts}\label{sec:basic}

We begin by recalling fundamental notions and definitions that will be used throughout the paper.

Let $S=\C[X_0,\dots,X_n]$ denote the standard graded polynomial ring over the field of complex numbers. We write $\dP^n$ for the projective space dual to $\P^n$, with polynomial ring $\check{S} := \C[Y_0,\dots,Y_n].$ Within $\P^n$, let $Z$ be a finite set of points. In all examples considered later, $Z$ will arise as a special configuration of points. To the subset $Z \subset \P^n$, we associate the dual hyperplane arrangement $\A_Z=\A=\{H_1,\ldots,H_{|Z|}\},$ where each $H_i \subset \dP^n$ is the hyperplane corresponding to the point $P_i \in Z$. If $P_i=(p_{i0}:\ldots:p_{in})$, then the hyperplane $\check{P_i}=H_i$ is defined by the linear equation
$$\ell_i :\quad p_{i0}Y_0+\cdots+p_{in}Y_n=0.$$

Let $\H\subset \P^n$ be a general linear subspace of codimension two. Its dual $\check\H\subset \dP^n$ is a line, which we denote by $\eL\subset \dP^n$. We adopt the following convention for coordinates:
\begin{itemize}
   \item points in $\P^n$ have homogeneous coordinates $X=(X_0:\ldots:X_n)$,
   \item points in $\dP^n$ have homogeneous coordinates $Y=(Y_0:\ldots:Y_n)$.
\end{itemize}
For notational convenience, we will use the shortened notation $P_i\circ Y=0$ to denote the defining equation of $H_i$, and similarly $A\circ X=0$, $B\circ X=0$ for the equations of $\H$ defined by points $A,B\in\dP^n$.

\medskip

Let $\Der(\check{S})$ denote the module of $\C$-linear derivations of $\check{S}$ in the variables $Y_0,\dots,Y_n$. The Euler derivation is defined as
$$\theta_E=\sum_{i=0}^n Y_i \frac{\partial}{\partial Y_i}\in\Der(\check{S}).$$
Given a homogeneous polynomial $f\in\check{S}$, the module of logarithmic derivations along $f$ is
$$\D(f)=\left\{\theta\in \Der(\check{S}) \;\middle|\; \theta(f)\in f\check{S} \right\}.$$
It is a classical fact that $\D(f) = \check{S}\theta_E \oplus \D_0(f),$
where $\D_0(f)$ is the kernel of the $\C$-linear map 
$$\Der(\check{S}) \to \check{S}, \qquad \theta \mapsto \theta(f).$$
The module $\D_0(f)$ is called the \emph{derivation module}, and the corresponding locally free sheaf is the \emph{derivation bundle} (see \cite{CHMN} for details). 
When $f$ is the defining polynomial of a hyperplane arrangement $\A_Z$, i.e.\ $f=\prod_{i=1}^{|Z|}\ell_i$, we will also write $\D(\A)$ in place of $\D(f)$.

\medskip

In the same spirit, we define $\Der^k(\check{S})$, the module of $k$-derivations, obtained as products of first-order differential operators with polynomial coefficients.

\begin{definition}\label{def:der-ka}
The module of $k$-derivations of $\check{S}$ is defined as
$$
   \Der^k(\check{S})
   =\left\{\, 
      \theta=\sum_{j=0}^{\binom{k+n}{n}-1} f_{\beta_j}(Y)\,\partial_{\beta_j} 
      \;\middle|\; f_{\beta_j}\in\check{S} 
   \,\right\},
$$
where for each $j$ the multi-index $\beta_j=(\beta_{j0},\dots,\beta_{jn})$ satisfies $|\beta_j|=k$. The operator $\partial_{\beta_j}$ is defined as the ordered product of first-order partial derivatives,
$$
   \partial_{\beta_j} 
   = \underbrace{\frac{\partial}{\partial Y_0}\cdot \ldots \cdot \frac{\partial}{\partial Y_0}}_{\beta_{j0}\ \text{times}}
     \underbrace{\frac{\partial}{\partial Y_1}\cdot \ldots \cdot \frac{\partial}{\partial Y_1}}_{\beta_{j1}\ \text{times}}
     \;\cdots\;
     \underbrace{\frac{\partial}{\partial Y_n}\cdot \ldots \cdot \frac{\partial}{\partial Y_n}}_{\beta_{jn}\ \text{times}}.
$$
\end{definition}

\medskip

For a homogeneous polynomial $f\in\check{S}$, we write $f_{Y_i}=\frac{\partial f}{\partial Y_i}$, and denote by
$$J=\langle f_{Y_0},\dots,f_{Y_n}\rangle$$
the \emph{Jacobian ideal} of $f$. Throughout, we abbreviate 
$$N = \binom{n+k}{n}-1,$$
and we assume that the coefficients $f_{\beta_j}$ are homogeneous polynomials of degree $d$.

\medskip


\section{Syzygies-based constructions}\label{sec:construction-generalized}

We now turn to the syzygy approach. This section recalls a classical construction and introduces its generalization. Our ultimate goal is to relate $k$-derivations to suitable modules of syzygies, extending results of \cite{CHMN,MT-G}, and \cite{Trok}.

Let $Z=\{P_1,\dots,P_s\}$ be a set of distinct points in $\P^2$, and let $P$ be a general point. Define $f=\check P_1\cdot\ldots \cdot \check P_s$, and let $J$ denote its Jacobian ideal. In \cite{MT-G}, a construction from \cite{CHMN} was generalized as follows. The original procedure uses degree-$d$ syzygies of $J$, restricted to the line $\check{P} \subset\check\P^2$, to produce a curve in $\P^2$ of degree $d+1$, passing through $Z$ and through $P$ with multiplicity $d$. The generalized construction instead employs degree-$d$ syzygies of $J^k$, yielding a curve of degree $d+k$, again passing through $Z$ and through $P$ with multiplicity $d$. In what follows we present a further generalization, which extends both of the above and also parts of the results of Trok \cite{Trok}.

\medskip

We fix the following notation. Let $Z=\{P_1,\dots,P_r\} \subset \P^n$ with
$$P_i = (x_{i0} : x_{i1} : \ldots : x_{in}),$$
and let $\check{P_i}$ denote the dual hyperplane. 
We write
$$f = \prod_{i=1}^{r} \check P_i$$
for the defining polynomial of $\A_Z$, and $J_f$ for its Jacobian ideal. We also write
$$X^{k}:=\bigl( X^\beta \;\big|\; \beta=(\beta_0,\dots,\beta_n)\in\N^{n+1},\; |\beta|=k \bigr),$$
where 
$$X^\beta = X_0^{\beta_0}X_1^{\beta_1}\cdots X_n^{\beta_n}.$$

Let $A,B\in\dP^n$ be generic points. Each of them defines a hyperplane in $\P^n$ via the equations $A\circ X=0$ and $B\circ X=0$, respectively.  Their intersection is a codimension-two 
subspace $\H\subset \P^n$, whose dual is a line $\eL=\check \H\subset \dP^n$, which will play a central role in the constructions below.

\begin{theorem}
\label{thm:syzConstr}
Let $Z=\{P_1,\dots,P_r\}$ be a finite set of points in $\P^n$, where 
$$P_i=(x_{i0}:\ldots:x_{in}), \qquad i=1,\dots,r,$$ and let $\A=\{\check P_1,\dots,\check P_r\}$
be the corresponding dual hyperplane arrangement in $\dP^n$.  
Fix two generic points $A, B\in\dP^n$, and let 
$$\eL = \{\, tA+(1-t)B \;\mid\; t\in\C \,\} \subset \dP^n$$
be the line they span. Consider a reduced syzygy $G=(g_{\beta_0},\dots,g_{\beta_N},g)$ of the ideal $J^k+( \eL )$, where $N=\binom{n+k}{n}-1$, each coefficient $g_{\beta_j}$ is homogeneous of degree $d$, and $Q\in\eL$ denotes a generic point of the line.  We associate to $G$ the polynomial $S_G(Q,X):= G(Q)\circ X^k,$ where $X^k=(X^\beta)_{|\beta|=k}$ denotes the vector of all monomials of degree $k$ in $X=(X_0,\dots,X_n)$.

Consider in $\dP^n$ the system of equations
\begin{equation}\label{s0}\tag{$\star$} 
   \begin{cases} 
      A\circ X=0,\\
      B\circ X=0,\\
      S_G(Q,X)=0.
   \end{cases}
\end{equation}
The system $(\star)$ will be called \emph{underdetermined} if its third equation reduces to
$$S_G(Q,X) = (Q\circ X)^k,$$
that is, if $S_G$ coincides with the $k$-th power of the linear form defined by $Q$.

Then the following holds:
\begin{enumerate}
   \item The system $(\star)$ is underdetermined only for points $Q$ lying on $\A \cap \eL$.
   \item All solutions of $(\star)$ lie on a hypersurface 
   $\mathcal{S}\subset\P^n$ of degree at most $d+k$.
   \item The hypersurface $\mathcal{S}$ contains the set $Z$.
   \item The hypersurface $\mathcal{S}$ has multiplicity at least $d$ along $\H$.
\end{enumerate}
\end{theorem}
\begin{proof}
(1) For a point $Q\in\dP^n$, let $J_Q$ denote the Jacobian vector 
$J=(f_{Y_0},\ldots,f_{Y_n})$ evaluated at $Q$. Suppose that the system \eqref{s0} is underdetermined for some $Q\in\eL$. By definition, this means
$$G(Q)\circ X^k=(Q\circ X)^k=0.$$
In particular, taking $X=J_Q$ we obtain
$$G(Q)\circ (J_Q)^k=(Q\circ J_Q)^k=0.$$
By Euler’s identity, $Q\circ J_Q=(\deg f)\, f(Q)$. Hence $(\deg f \cdot f(Q))^k=0$, which forces $f(Q)=0$. Therefore the system can be underdetermined only for points $Q\in \A\cap\eL$, 
as required.

(2)  
Let $Q\in\eL$. By definition, $\check{Q}$ contains $\H$, and it can be written as 
$$\check{Q} :\quad (t_QA+(1-t_Q)B)\circ X =0$$
for some $t_Q\in\C$. The parameter $t_Q$ can be expressed explicitly as
$$t_Q = \frac{B\circ X}{(B-A)\circ X}.$$
Now consider the polynomial $S_G(X)$. By definition,
\begin{align*}
   S_G(X)=G(Q)\circ X^k=G(t_QA+(1-t_Q)B)\circ X^k=G\bigl(t_Q(A-B)+B\bigr)\circ X^k.
\end{align*}
Substituting the value of $t_Q$ we obtain
\begin{align*}
   S_G(X)
   &=G\!\left( \left(\frac{B\circ X}{(B-A)\circ X}\right)(A-B)+B\right)\circ X^k \\
   &=G\bigl((B\circ X)(A-B)+((B-A)\circ X)\,B\bigr)\circ X^k \\
   &=G\bigl((B\circ X)\,A-(A\circ X)\,B\bigr)\circ X^k.
\end{align*}
By construction, each coefficient of the syzygy $G$ is homogeneous of degree $d$, while $X^k$ denotes the collection of all monomials of degree $k$. Hence $S_G(X)$ is a homogeneous polynomial in the coordinates of $X$ of degree at most $d+k$. Consequently, the vanishing locus of $S_G(X)$ is contained in a hypersurface $\mathcal{S}\subset \P^n$ of degree $\leq d+k$, and by definition all solutions of \eqref{s0} lie on $\mathcal{S}$. This proves the second claim.

(3)  
We show that each point $P_i\in Z$ lies on the hypersurface $\mathcal{S}$.  

Let $f=\check{P}_1\cdots \check{P}_r$ be the defining polynomial of the arrangement $\A=\{\check P_1,\dots,\check P_r\}$. For each $i$ let $Q_i:=\eL\cap \check{P}_i$. By construction, $\check P_i(Q_i)=0$. Recall that if $P_i=(x_{i0}:\ldots:x_{in})$, then $\check P_i$ is given by
$$\check P_i : \quad x_{i0}Y_0+\cdots+x_{in}Y_n=0.$$
Differentiating $f$ we obtain
$$\frac{\partial f}{\partial Y_j}
   =\sum_{i=1}^r x_{ij}\, \check P_1\cdots \widehat{\check {P_i}}\cdots \check P_r,
$$ where the hat indicates omission of the $i$-th factor. Evaluating at $Q_i$ gives
$$\frac{\partial f}{\partial Y_j}(Q_i)=x_{ij}\,\prod_{m\neq i}\check P_m(Q_i).$$
Therefore,
$$J(Q_i)=(x_{i0},\dots,x_{in})\cdot\prod_{m\neq i}\check P_m(Q_i)
   =P_i \cdot\prod_{m\neq i}\check P_m(Q_i).
$$
Since $\H$ is generic, the product $\prod_{m\neq i}\check P_m(Q_i)$ is nonzero. 
Hence
$$J^k(Q_i)=(P_i)^k \cdot \Bigl(\prod_{m\neq i}\check P_m(Q_i)\Bigr)^k.$$
Now $G$ is, by definition, a syzygy of $J^k$, so in particular
$G(Q_i)\circ P_i^k=0.$ This shows that $P_i\in\mathcal{S}$, and hence $\mathcal{S}$ passes through all points of $Z$.

(4) Let $U=(U_0,\dots,U_N)$ denote homogeneous coordinates on $\P^N$, where $N=\binom{n+k}{n}-1$. Fix a point $P\in \H$. The equation
$$P^k \circ U=0$$
defines a hyperplane in $\P^N$. For each $Q\in\eL$, the evaluation $G(Q)$ gives a point of $\P^N$, and as $Q$ varies along $\eL$ this defines a curve of degree $d$ in $\P^N$. By Bézout’s theorem, this curve intersects the hyperplane $P^k\circ U=0$ in exactly $d$ points (counted with multiplicities). Denote these intersection points by $u_1,\dots,u_d$, with $u_j=G(Q_j)$ for some $Q_j\in\eL$.  
By construction, we then have
$$G(Q_j)\circ P^k=0 \qquad \text{for } j=1,\dots,d.$$
This shows that the hypersurface $\mathcal S$ passes through $P$ with multiplicity at least $d$.
\end{proof}

\section{Derivations}\label{sec:derivations}

Recall Definition \ref{def:der-ka}, where we introduced the set of $k$-derivations $\Der^k(\check{S})$. We now define a subset of this set, associated with the dual hyperplane arrangement $\A_Z$.
\begin{definition}\label{def:de-ka}
We define the module of $k$-derivations associated with the arrangement $\A_Z$ as
$$
   D^k(\A_Z):= \{\, \theta \in \Der^k(\check{S}) \mid \theta(f) \in (J^{k-1}f) \,\},
$$
where $f$ denotes the defining equation of $\A_Z$.
\end{definition}
The goal of this section is to establish a correspondence between elements of $D^k(\A_Z)$ and certain polynomial functions on $\P^n$. Most of the material presented here is based on the (currently unpublished) work of Trok \cite{Trok}. We follow his line of reasoning, extending the argument from the case $k=1$ to arbitrary $k \geq 1$.

Let us recall the general setting. In $\P^n$ we consider the set of points $Z=\{P_1,\dots,P_r\}$ together with a generic codimension two subspace $\H$. Dually, in $\dP^n$ this corresponds to the set of hyperplanes $H_i=\check P_i$ for $i=1,\dots,r$, and the line $\eL=\check \H$. Our aim is to construct a transformation that associates to a $k$-derivation $\theta\in D^k(\A_Z)$, with polynomial coefficients $f_{\beta_j}(Y)$ of degree $d$, a map $F(X)$ of degree $d+k$. We begin by introducing some notation and recalling a few basic facts, following \cite[Section 4.3]{Trok}.

Let $b_0,\dots,b_{n-1}$ be $n$ vectors in $\dP^n$ spanning a hyperplane $H$, where
$$b_j = [A_{j0},A_{j1},\dots,A_{jn}], \qquad j=0,\dots,n-1.$$
For each $i=0,\dots,n$, define
$$\mu_i:=\det\!\begin{bmatrix}
      0 & \cdots & 1 & \cdots & 0 \\[2pt]
      A_{00} & \cdots & A_{0i} & \cdots & A_{0n} \\
      A_{10} & \cdots & A_{1i} & \cdots & A_{1n} \\
      \vdots & \ddots & \vdots & \ddots & \vdots \\
      A_{n-1,0} & \cdots & A_{n-1,i} & \cdots & A_{n-1,n}
   \end{bmatrix}.$$
Now let $M_i$ denote the determinant obtained from $\mu_i$ by replacing the coordinates of $b_0$ with the variables $X=(X_0,\dots,X_n)$, that is,
$$
   M_i:=\det\!\begin{bmatrix}
      0 & \cdots & 1 & \cdots & 0 \\[2pt]
      X_0 & \cdots & X_i & \cdots & X_n \\
      A_{10} & \cdots & A_{1i} & \cdots & A_{1n} \\
      \vdots & \ddots & \vdots & \ddots & \vdots \\
      A_{n-1,0} & \cdots & A_{n-1,i} & \cdots & A_{n-1,n}
   \end{bmatrix}.
$$
Then, we set
$$M=(M_0:\ldots :M_n).$$ With this notation in place, we can now define the map
\begin{equation}\label{eq:def eta}
   \eta : \theta=\sum_{j=0}^{N} f_{\beta_j}(Y)\,\partial_{\beta_j}
      \;\longmapsto\;
      F=\sum_{j=0}^{N} f_{\beta_j}(M)\,X^{\beta_j},
\end{equation}
where $X^{\beta_j}=X_0^{\beta_{j0}}\cdots X_n^{\beta_{jn}}$. 

Note that for any point $Q=(q_0:\ldots:q_n)\in \P^n$, with dual hyperplane
$$\check{Q}: \; \sum_{i=0}^n q_i Y_i=0,$$
we have $\sum_{i=0}^n q_i M_i=0$ at $Q$. 
This follows directly from the definition of $M_i$ and the properties of determinants. Furthermore, observe that the function $F=\eta(\theta)$ may be regarded as depending on the variables $X$ and on the coefficients $A_{ij}$, for $1\leq i\leq n-1$ and $0\leq j\leq n$. By substituting the vectors $b_i$ (for $i\geq 1$) with points of $\P^n$, we obtain a function defined on $X\in \P^n$. More precisely, let 
$$
   \alpha = \{\alpha_1,\dots,\alpha_{n-1}\}, \qquad 
   \alpha_i = (\alpha_{i0}:\ldots:\alpha_{in}),
$$
be a set of linearly independent vectors in $\P^n$. 
For each $i=0,\dots,n$, define
$$
   e_{\alpha}(M_i):= 
   \det\!\begin{bmatrix}
      0 & \cdots & 1 & \cdots & 0 \\[2pt]
      X_0 & \cdots & X_i & \cdots & X_n \\
      \alpha_{10} & \cdots & \alpha_{1i} & \cdots & \alpha_{1n} \\
      \vdots & \ddots & \vdots & \ddots & \vdots \\
      \alpha_{n-1,0} & \cdots & \alpha_{n-1,i} & \cdots & \alpha_{n-1,n}
   \end{bmatrix}.
$$
We then write 
$$e_{\alpha}(M):=\bigl(e_{\alpha}(M_0):\ldots:e_{\alpha}(M_n)\bigr),$$
and define
$$e_{\alpha}(F):= \sum_{j=0}^{N} f_{\beta_j}\bigl(e_{\alpha}(M)\bigr)\, X^{\beta_j}.$$
Finally, we denote by $\H:=\mathrm{span}(\alpha)$ the subspace spanned by the vectors in $\alpha$; the choice of the symbol $\H$ here is intentional.

\medskip
Lemma \ref{le:trok4.9} below is a generalization of \cite[Lemma 4.9]{Trok}; it will be used in the sequel of this section.

\begin{lemma}\label{le:trok4.9}
With the above notation, for every $Q \in\dP^n$ such that $\check{Q}$ vanishes on $\H$, there exists a nonzero linear form $h$ vanishing on $\H$ such that, modulo $\check{Q}$,
$$e_{\alpha}(F)\;=\; h(X)^d \sum_{j=0}^{N} f_{\beta_j}(Q)\, X^{\beta_j}.$$
\end{lemma}

\begin{proof}
We may change coordinates in $\P^n$ so that
\[
   \alpha_1=(0:0:1:0:\ldots:0),\quad 
   \alpha_2=(0:0:0:1:0:\ldots:0),\quad\dots,\quad
   \alpha_{n-1}=(0:\ldots:0:1),
\]
and $\check{Q}$ is given by $X_0=0$, hence $Q=(1:0:\dots:0)$. 
It is straightforward to check that
$$
   e_{\alpha}(M_0)=X_1,\qquad 
   e_{\alpha}(M_1)=-X_0,\qquad 
   e_{\alpha}(M_j)=0\quad\text{for } j>1.
$$
Therefore,
$$e_{\alpha}(F)=\sum_{j=0}^{N}f_{\beta_j}(e_{\alpha}(M))\,X^{\beta_j}=\sum_{j=0}^{N} f_{\beta_j}(X_1: -X_0:0: \ldots :0)\,X^{\beta_j}.$$
Since $\check{Q}$ imposes $X_0=0$, modulo $\check{Q}$ we obtain
$$e_{\alpha}(F)=\sum_{j=0}^{N} f_{\beta_j}(X_1:0:0:\ldots :0)\, X^{\beta_j}=X_1^d \sum_{j=0}^{N} f_{\beta_j}(1:0:0:\ldots :0)\,X^{\beta_j}.
$$
Thus, taking $h(X)=X_1$ and recalling that $Q=(1:0:\ldots:0)$, the claim follows.
\end{proof}

\begin{remark}\label{re:vanishing on span alpha}
Lemma~\ref{le:trok4.9} implies that $e_{\alpha}(F)$ vanishes on $\H$ with multiplicity at least $d$.
\end{remark}

The following result is stated in \cite[Proposition 4.12]{Trok}.
\begin{proposition}\label{pro:C-algebry}
Let $\H=\mathrm{span}(\alpha)$ and $\eL=\check\H$, as above. 
Then there is an isomorphism of $\C$-algebras
$$\C[Y]/I(\eL) \;\cong\; \sym_{\C}[I(\H)].$$
\end{proposition}
\begin{proof}
Using the change of variables from the proof of Lemma \ref{le:trok4.9}, 
it suffices to compare the degree one components:
$$\bigl(\C[Y]/(Y_2,\dots,Y_n)\bigr)_1 \;\cong\; \sym_{\C}[X_0,X_1]_1.$$
Since both algebras are generated in degree $1$, the claim follows.
\end{proof}

The next definition is a $k\geq 1$ generalization of \cite[Definition 4.13]{Trok}.
\begin{definition}\label{def:trok4.13module}
Let $I(Z)$ denote the ideal of a set $Z\subset\P^n$, and let 
$$\mathcal{P}:=\sym_{\C}[X_0,X_1].$$
We define a graded $\mathcal{P}$-module $I(Z,\H)\subset \C[X]$ by prescribing its 
$(d+k)$-th graded component as
$$\bigl[I(Z)\cap I(\H)^d\bigr]_{d+k}.$$
\end{definition}

We now need to restrict the differentials $D^k(\A_Z)$ to the line $\eL$. In degree $d$, this restriction can be expressed as
$$H^0\bigl(\,\widetilde{D^k(\A_Z)}\otimes \mathcal{O}_{\eL}(-d),\, \eL \bigr),$$ 
but this perspective is not useful for our purposes. Instead, we require a generalization of Trok’s Proposition 3.13.

\begin{proposition}\label{pro:trok3.13restr-diff}
Let $f=0$ be the defining equation of $\A_Z$, written as $f = \check P_1 \cdot \ldots \cdot \check P_r.$
Then
$$
   D^k(\A_Z)=
   \left\{\, \theta \in \Der^k\bigl(\C[Y], \C[Y]/I(\eL)\bigr) 
      \;\middle|\; \theta(\check P_i) \in (\overline{\check P_i}), \; i=1,\dots,r \,\right\},
$$
where $\overline{\check P_i}$ denotes the class of $\check P_i$ in $\C[Y]/I(\eL)$.
\end{proposition}

\begin{proof}
The argument follows the same line as in \cite{Trok}; we include it here for completeness.

\smallskip
\noindent
\emph{Case $r=1$.}  
Assume $H_1=H=\{Y_0=0\}$. 
Then $D^k(H)$ consists of all $k$-derivations $\theta$ such that
$$
   \theta(Y_0) 
     =\sum_{j=1}^N f_{\beta_j}(Y)\,\partial_{\beta_j}(Y_0)\in (Y_0),
$$
(with $J=1$ in this case). 
Hence $D^k(H)$ is generated by $Y_0\partial_{(k,0,\dots,0)}$ together with all remaining $\partial_{\beta_j}$. Restricting this module to $\eL$ yields a $\C[Y]/I(\eL)$-module with basis
$$\bigl\{ (Y_0+I(\eL))\partial_{(k,0,\dots,0)}, \; \partial_{\beta_j} \ (j>0) \bigr\}.$$
This proves the claim when $\A_Z=H$.

\smallskip
\noindent
\emph{Case $r>1$.}  
For each $i=1,\dots,r$, set $U_i := \A_Z \setminus H_i$. 
By the genericity of $\eL$, the open sets $\{U_i\}_{i=1}^r$ form an open cover of $\eL$. 

Let $\sigma \in H^0\bigl(\widetilde{D^k(\A_Z)}\otimes \mathcal{O}_{\eL}(-d), \eL\bigr)$. 
Then
$$
   \sigma|_{U_i \cap\eL} \in H^0\bigl(\widetilde{D^k(\A_Z)}\otimes\mathcal{O}_{\eL}(-d), \eL \cap U_i\bigr).
$$
Moreover, on each such open subset we have
$$
   H^0\bigl(\widetilde{D^k(\A_Z)}\otimes\mathcal{O}_{\eL}(-d), \eL\cap U_i\bigr)
   =H^0\bigl(\widetilde{D^k(H_i)}\otimes \mathcal{O}_{\eL}(-d), \eL\cap U_i\bigr).
$$
For $i\neq j$ one further obtains
$$
   H^0\bigl(\widetilde{D^k(H_i)}\otimes \mathcal{O}_{\eL}(-d), \eL \cap U_j\bigr)
   =H^0\bigl(\Der^k \otimes \mathcal{O}_{\eL}(-d), U_j \cap \eL\bigr).
$$

Thus the following conditions are equivalent:
\begin{enumerate}
   \item $\sigma \in H^0\bigl(\widetilde{D^k(\A_Z)}\otimes \mathcal{O}_{\eL}(-d), \eL\bigr)$,
   \item $\sigma|_{U_i \cap \eL} \in H^0\bigl(\widetilde{D^k(H_i)}\otimes \mathcal{O}_{\eL}(-d), \eL \cap U_i\bigr)$ for all $i=1,\dots,r$,
   \item $\sigma \in H^0\bigl(\widetilde{D^k(H_i)}\otimes \mathcal{O}_{\eL}(-d), \eL\bigr)$ for all $i=1,\dots,r$.
\end{enumerate}
From the first part of the proof, condition (3) is equivalent to
\begin{enumerate}
   \item[4.] $\sigma \in 
   \bigl\{\theta \in \Der^k(\C[Y], \C[Y]/I(\eL)) 
      \;\big|\; \theta(\check P_i) \in (\overline{\check P_i}),\ i=1,\dots,r \bigr\}$.
\end{enumerate}
This completes the proof.
\end{proof}

\begin{corollary}
Thus, $D^k(\A_Z)|_{\eL}$ is a $\C[Y]/I(\eL)$-module. After tensoring with $\mathcal{P}$, it may therefore be regarded as an $\mathcal{P}$-module.  
\end{corollary}

Recall that $\theta_E$ denotes the Euler derivation
$$\theta_E=\sum_{i=0}^n Y_i \frac{\partial}{\partial Y_i}.$$
In the case $k=1$, this derivation generates the kernel of the map defined in~\eqref{eq:def eta}. 
We now generalize this construction to arbitrary $k\geq 1$. Recall that a multi-index is written as 
$$\beta_j = (\beta_{j0},\dots,\beta_{jn}), \qquad |\beta_j| = k,$$
with 
$$
   \partial_{\beta_j}=
      \Bigl(\frac{\partial}{\partial Y_0}\Bigr)^{\beta_{j0}}
      \Bigl(\frac{\partial}{\partial Y_1}\Bigr)^{\beta_{j1}}
      \cdots
      \Bigl(\frac{\partial}{\partial Y_n}\Bigr)^{\beta_{jn}},
$$
see Definition \ref{def:der-ka}. Now take a multi-index $\gamma_t=(\gamma_{t0},\dots,\gamma_{tn})$ with $|\gamma_t|=k-1$. 
For each $m=0,\dots,n$, define 
$$\beta_{j_m}:=\gamma_t+e_m,$$
where $e_m$ denotes the standard unit vector of length $n+1$ with $1$ in position $m$ and $0$ elsewhere.

\begin{definition}\label{de:Eulery}
Let $\gamma_t$ and the indices $\beta_{j_m}$ be as above.  
We define the $k$-Euler derivation associated with $\gamma_t$ as
$$\theta_{E_{\gamma_t}} := \sum_{m=0}^n Y_m \, \partial_{\gamma_t + e_m}.$$
\end{definition}
\begin{remark}
Informally, $\theta_{E_{\gamma_t}}$ may be viewed as an analogue of the Euler derivation $\theta_E$, multiplied by a common factor of derivations $\partial_{\gamma_t}$. Such derivations always occur among the minimal generators of the syzygy module, but they do not carry essential information. This motivates the following definition, which is the natural analogue of $D_0(\A_Z)$ in the case $k=1$.
\end{remark}
\begin{definition}\label{de:D0}
We define
\[
   D_0^k(\A_Z) := D^k(\A_Z)\,\Big/\, 
   \bigoplus_{t=1}^{\binom{n+k-1}{n}} \check{S} \,\theta_{E_{\gamma_t}}.
\]
\end{definition}
\begin{example}
Let $n=2$ and $k=2$. 
For the multi-index $\gamma=(1,0,0)$ we obtain
$$
   \theta_{E_{\gamma}} 
   =Y_0 \,\partial_{(2,0,0)}+Y_1 \,\partial_{(1,1,0)}+Y_2 \,\partial_{(1,0,1)}.
$$
This illustrates the general construction: $\theta_{E_{\gamma}}$ is obtained by multiplying the classical Euler derivation $\theta_E=Y_0\partial_{(1,0,0)}+Y_1\partial_{(0,1,0)}+Y_2\partial_{(0,0,1)}$ by the common differential factor $\partial_{\gamma}=\partial_{(1,0,0)}$.
\end{example}
\begin{remark}\label{re:trok3.13restr-diff-part2}
The statement of Proposition \ref{pro:trok3.13restr-diff} remains valid if $D^k(\A_Z)$ is replaced by $D_0^k(\A_Z)$.
\end{remark}
\begin{proof}
This follows from the general fact that if $B$ is a submodule of a module $C$, then
$$(C/B)\otimes {\cal O}_{\dP^n}(-d) \;\cong\; 
   (C\otimes {\cal O}_{\dP^n}(-d)) / (B\otimes {\cal O}_{\dP^n}(-d)).$$
Applying this with $C=D^k(\A_Z)$ and $B=\bigoplus_{t} \check S \,\theta_{E_{\gamma_t}}$ gives the claim. The remainder of the argument is identical to the proof of Proposition~\ref{pro:trok3.13restr-diff}.
\end{proof}
Now we may formulate the main result of this section, the generalization of Theorems 4.8 and 4.14 from \cite{Trok}.
\begin{theorem}\label{the:trok4.8and4.14}
The map
$$
   \eta :\;
   \theta=\sum_{j=0}^{N} f_{\beta_j}(Y)\,\partial_{\beta_j}
      \;\longmapsto\;
      e_{\alpha}(F) = \sum_{j=0}^{N} f_{\beta_j}\bigl(e_{\alpha}(M)\bigr)\,X^{\beta_j},
$$
defined as in~\eqref{eq:def eta} with substitution $e_{\alpha}$ for a general basis $\alpha$ of $\H$, is a morphism of graded $\mathcal{P}$-modules. 
Its kernel is 
$$\bigoplus_{t=1}^{\binom{n+k-1}{n}} \check S \,\theta_{E_{\gamma_t}},$$
so that there is an isomorphism
$$D_0^k(\A_Z)|_{\eL} \;\cong\; I(Z,\H)(-k).$$
\end{theorem}
\begin{proof}
We first reduce to the case $r=1$. Indeed, if $Z=\{P_1,\dots,P_r\}$ and $H_i=\check P_i$, then
$$D_0^k(H_i)|_{\eL} \;\cong\; I(P_i,\H)(-k), \qquad i=1,\dots,r,$$
and the general case follows by intersecting over all $i$.

\smallskip
\noindent
\emph{Step 1. The image of $\eta$.}  
Let $Z=\{P\}$ with $P=(p_0,\dots,p_n)$ and let 
$$H=\{\psi(Y)=p_0Y_0+\cdots+p_nY_n= P\circ Y =0 \}.$$
For $\theta=\sum f_{\beta_j}(Y)\,\partial_{\beta_j}$ we have
$$
   \theta \in D^k(\A_Z)|_{\eL}
   \;\Longleftrightarrow\;
   \sum_{j=0}^{N} f_{\beta_j}(Y)\,\partial_{\beta_j}(\psi)\in (\psi).
$$
Since $\partial_{\beta_j}(\psi)=P^{\beta_j}$, this is equivalent to
$$\sum_{j=0}^{N} f_{\beta_j}(Y)\,P^{\beta_j} \;\in\; (\psi).$$
Evaluating at $Y=Q\in \eL$, and applying Lemma~\ref{le:trok4.9}, we obtain
$$h(Q)^d \sum_{j=0}^{N} f_{\beta_j}(Q)\,P^{\beta_j}=0,$$
where $h$ is a nonzero linear form vanishing along $\H$. 
Hence $e_{\alpha}(F)$ of degree $d+k$ vanishes at $P$ and with multiplicity $d$ along $\H$, so
$$\eta(\theta)=e_{\alpha}(F) \in I(P,\H).$$

\smallskip
\noindent
\emph{Step 2. The kernel of $\eta$.}  
Suppose $e_{\alpha}(F)=0$ for a generic basis $\alpha$ of $\H$. 
Then for $Q\in \eL$,
$$h(X)^d \sum_{j=0}^{N} f_{\beta_j}(Q)\,X^{\beta_j}=0 \pmod{\check{Q}},$$
which implies
\begin{equation}\label{eq:kernel}
   \sum_{j=0}^{N} f_{\beta_j}(Q)\,X^{\beta_j}\;=\; (Q\circ X)\cdot W(X),
\end{equation}
for some homogeneous polynomial $W$ of degree $k-1$. Choose a monomial $X^{\gamma_t}$ of degree $k-1$. 
As noted before, there exist indices $\beta_{j_m}$, $m=0,\dots,n$, with
$$\beta_{j_m} = \gamma_t + e_m,$$
where $e_m$ is the $m$-th standard basis vector. 
Then the left-hand side of~\eqref{eq:kernel} may be written as
$$\sum_{j=0}^{N} f_{\beta_j}(Q)\,X^{\beta_j}
   =X^{\gamma_t}\sum_{m=0}^{n} f_{\beta_{j_m}}(Q)\,X_m \;+\; V(X),
$$
with $V(X)$ not divisible by $X^{\gamma_t}$. 
Similarly, the right-hand side expands as
$$(Q\circ X)\cdot W(X) 
   =c\,X^{\gamma_t}(q_0X_0+\cdots+q_nX_n) \;+\; U(X),
$$
where $Q=(q_0:\ldots:q_n)$ and $U(X)$ is not divisible by $X^{\gamma_t}$. Comparing coefficients of $X^{\gamma_t}X_m$ on both sides yields
$$f_{\beta_{j_m}}(Q) = c q_m.$$
Thus the map
$$Q \;\longmapsto\; (f_{\beta_{j_0}}(Q):\cdots:f_{\beta_{j_n}}(Q))$$
is projectively the identity. Repeating this argument for all $\gamma_t$ shows that $e_{\alpha}(F)$ arises from a linear combination of the Euler-type derivations $\theta_{E_{\gamma_t}}$. Hence the kernel of $\eta$ is precisely
$$\bigoplus_{t=1}^{\binom{n+k-1}{n}} \check S \,\theta_{E_{\gamma_t}},$$
which completes the proof.
\end{proof}

\begin{corollary}\label{co:dimension}
    From Theorem \ref{the:trok4.8and4.14} we obtain that
 $$\dim(I(Z,\H)_{d+k})=\sum_{i=1}^{\binom{k+n}{n}} \max\{0,(d+k-a_i)\},$$
 where  $D_0^k(\A_Z)|_{\eL}$ splits as $\bigoplus_{i=0}^{\binom{k+n}{n}}{\cal O}(-a_i)$.
\end{corollary}

\begin{remark}\label{re: syz-D}

There is an isomorphism between $\syz J^k$ and $D^k(\A_Z)$, as is the case between $\syz J$ and $D(\A_Z)$.
\end{remark}
Indeed, it is enough to consider the following diagram, analogous to the diagram for $n=2$, see \cite{MT-G}.

\begin{center}
\footnotesize
	\begin{tikzpicture}
		\matrix (m) [matrix of math nodes,row sep=2.5em,column sep=3.5em,minimum width=2em,nodes={minimum height=3em, anchor=center}]
			{
					&  & 0 & 0 & \\
					& 0 & \syz((J(|Z|-1))^k) & D^k(\A) & \\
					0 & J^{k-1}f\check S(-1) & J^{k-1}f\check S(-1)\oplus \check S^{\binom{k+n}{n}} & \check S^{\binom{k+n}{n}} & 0\\
					0 & J^{k-1}f\check S(-1) & (J(|Z|-1) )^k& ((J/f\check S) (|Z|-1))^k & 0\\
					& 0 & 0 & 0 & \\};
			\path[-stealth]
				(m-3-1) edge node [left] {} (m-3-2)
				(m-3-2) edge node [left] {} (m-3-3)
				(m-3-3) edge node [left] {} (m-3-4)
				(m-3-4) edge node [left] {} (m-3-5)
				(m-4-1) edge node [left] {} (m-4-2)
				(m-4-2) edge node [above] {$\times f$} (m-4-3)
				(m-4-3) edge node [left] {} (m-4-4)
				(m-4-4) edge node [left] {} (m-4-5)
				(m-2-2) edge node [left] {} (m-3-2)
				(m-3-2) edge node [left] {} (m-4-2)
				(m-4-2) edge node [left] {} (m-5-2)
				(m-1-3) edge node [left] {} (m-2-3)
				(m-2-3) edge node [left] {} (m-3-3)
				(m-3-3) edge node [left] {} (m-4-3)
				(m-4-3) edge node [left] {} (m-5-3)
				(m-1-4) edge node [left] {} (m-2-4)
				(m-2-4) edge node [left] {} (m-3-4)
				(m-3-4) edge node [right] {} (m-4-4)
				(m-4-4) edge node [left] {} (m-5-4)
				;
	\end{tikzpicture}.
\end{center}

\section{Unexpectedness}\label{sec:unexpectedness}

In this section we prove sufficient conditions ensuring that the hypersurface constructed in Section~\ref{sec:construction-generalized} is unexpected. The result is the natural analogue of Theorem 5.5 (concerning unexpected curves) from \cite{MT-G}.

\begin{theorem}\label{the:unexp}
Let $Z\subset\dP^n$ and $\H\subset\P^n$, together with the dual line $\eL=\check\H\subset\dP^n$, and let $f$ be as above. Suppose the splitting type of the $k$-th derivation bundle is
$$
(a_1,a_2,\dots,a_N)=
(a,\dots,a,\,
 a+\epsilon_1,\dots,a+\epsilon_1,\,
 a+\epsilon_2,\dots,a+\epsilon_2,\,\dots\, ,
 a+\epsilon_s),
$$
where $\epsilon_0=0$, $1\leq \epsilon_1<\epsilon_2<\cdots<\epsilon_s$, and each $a+\epsilon_i$ appears with multiplicity $t_i$, so that $t_0+\cdots+t_s=N$. Then the hypersurface $\Xi$ of type $(a+\epsilon_j+k,\,a+\epsilon_j)$ is unexpected if
\begin{equation}\label{eq:teza 5.1}
   0<\sum_{i=j+1}^s t_i(\epsilon_i-\epsilon_j-1).
\end{equation}
\end{theorem}
    
\begin{proof}
Let us write $I(Z,\H)=I(Z)\cap I(\H)^{a+\epsilon_j}$. By Corollary \ref{co:dimension} we have
\begin{align*}
   \dim[I(Z,\H)]_{a+\epsilon_j+k}
   &=\sum_{\ell=1}^{\binom{k+n}{n}} \max\{0,\,a+\epsilon_j-a_\ell+1\} \\
   &=(\epsilon_j+1)t_0+(\epsilon_j+1-\epsilon_1)t_1+\cdots+(\epsilon_j+1-\epsilon_j)t_j.
\end{align*}
On the other hand, the expected dimension is
$$
   \binom{a+\epsilon_j+k+2}{2}
   -|Z|
   -\sum_{i=0}^{d-1}(i+1)\binom{d+k-i+n-2}{n-2},
$$
where the last sum counts the number of conditions imposed on forms of degree $d+k$ by vanishing along a codimension-two subspace of $\P^n$ with multiplicity $d$ (see \cite[Lemma~2.1]{CHMN}). Thus, to prove the theorem it suffices to show that condition \eqref{eq:teza 5.1} is equivalent to the inequality
$$
   \dim[I(Z,\H)]_{a+\epsilon_j+k}>
   \binom{a+\epsilon_j+k+2}{2}
   -|Z|
   -\sum_{i=0}^{d-1}(i+1)\binom{d+k-i+n-2}{n-2}.
$$

\begin{lemma}\label{le:compact-formula}
For integers $d,k,n\geq 1$ we have
\begin{align*}
\sum_{i=0}^{d-1}(i+1)\binom{d+k-i+n-2}{n-2}
&= (d+k+n)\left[\binom{d+k+n-1}{n-1}-\binom{k+n-1}{n-1}\right] \\
&\quad -(n-1)\left[\binom{d+k+n}{n}-\binom{k+n}{n}\right].
\end{align*}
\end{lemma}

\begin{proof}
Set
$$S=S(d,k,n)=\sum_{i=0}^{d-1}(i+1)\binom{d+k-i+n-2}{n-2}.$$
Put $\iota=i+1$, so
$$S=\sum_{\iota=1}^d \iota\binom{d+k+n-1-\iota}{n-2}.$$
Now set $t=d+k+n-1-\iota$; then
$$S=\sum_{t=d+k+n-1-d}^{d+k+n-2}(d+k+n-1-t)\binom{t}{n-2}.$$
Let $A=d+k+n-1$. Then
$$S=
A\sum_{t=A-d}^{A-1}\binom{t}{n-2}-\sum_{t=A-d}^{A-1}t\binom{t}{n-2}.$$
Using standard binomial identities:
$$
\sum_{u=0}^{M}\binom{u}{r}=\binom{M+1}{r+1},\qquad
\sum_{u=0}^{M}u\binom{u}{r}=(r+1)\binom{M+2}{r+2}-\binom{M+1}{r+1},
$$
we obtain
$$ A\sum_{t=A-d}^{A-1}\binom{t}{n-2}-\sum_{t=A-d}^{A-1}t\binom{t}{n-2}=
$$
$$=A\sum_{t=0}^{A-1}\binom{t}{n-2}-A\sum_{t=0}^{A-d-1}\binom{t}{n-2}-\sum_{t=0}^{A-1}t\binom{t}{n-2}+\sum_{t=0}^{A-d-1}t\binom{t}{n-2}=
$$
$$=A\binom{A}{n-1}-A\binom{A-d}{n-1}-(n-1)\binom{A+1}{n}+\binom{A}{n-1}+(n-1)\binom{A-d+1}{n}-\binom{A-d}{n-1.}$$
Thus,
$$
S=(A+1)\left[\binom{A}{n-1}-\binom{A-d}{n-1}\right]
-(n-1)\left[\binom{A+1}{n}-\binom{A-d+1}{n}\right].
$$
Substituting $A=d+k+n-1$ yields the desired formula.
\end{proof}

From the formula for $c_1(D_0^k(\A_Z))$ we know
$$|Z|=\binom{k+n-1}{n}+\sum_{i=0}^s (a+\epsilon_i)t_i.$$
Hence, we must analyse the inequality
\begin{align*}
\dim[I(Z,\H)]_{a+\epsilon_j+k}>\;&
\binom{a+\epsilon_j+k+2}{2}-|Z|\\
&-\Bigl((a+\epsilon_j+k+n)\Bigl[\binom{a+\epsilon_j+k+n-1}{n-1}-\binom{k+n-1}{n-1}\Bigr] \\ &
-(n-1)\Bigl[\binom{a+\epsilon_j+k+n}{n}-\binom{k+n}{n}\Bigr]\Bigr).
\end{align*}
For brevity, set $d=a+\epsilon_j$. 
Consider
\begin{align*}
\textbf{P}&:=\binom{d+k+n}{n}-S=\binom{d+k+n}{n}\\
&\quad-\left[(d+k+n)\left(\binom{d+k+n-1}{n-1}-\binom{k+n-1}{n-1}\right)
-(n-1)\left(\binom{d+k+n}{n}-\binom{k+n}{n}\right)\right].
\end{align*}
Grouping the terms with $\binom{d+k+n}{n}$ and using $n\binom{N}{n}=N\binom{N-1}{n-1}$ with $N=d+k+n$, we obtain
$$\textbf{P}=(d+k+n)\binom{k+n-1}{n-1}-(n-1)\binom{k+n}{n}.$$
Subtracting $|Z|$ gives
$$\textbf{P}-|Z|=(d+k+n)\binom{k+n-1}{n-1}-(n-1)\binom{k+n}{n}
-\binom{k+n-1}{n}-\sum_{i=0}^s (a+\epsilon_i)t_i.
$$
Now using the identities
$$
   \binom{k+n-1}{n}=\frac{k}{n}\binom{k+n-1}{n-1}, \qquad
   \binom{k+n}{n}=\frac{k+n}{n}\binom{k+n-1}{n-1},
$$
we simplify to
$$\textbf{P}-|Z|=(d+1)\binom{k+n-1}{n-1}-\sum_{i=0}^s (a+\epsilon_i)t_i.$$
Substituting back $d=a+\epsilon_j$ gives
$$\textbf{P}-|Z|=(a+\epsilon_j+1)\binom{k+n-1}{n-1}-\sum_{i=0}^s (a+\epsilon_i)t_i.$$
Finally, observe that
$$
(a+\epsilon_j+1)\binom{k+n-1}{n-1}-\sum_{i=0}^s (a+\epsilon_i)t_i<\sum_{i=0}^j (\epsilon_j+1-\epsilon_i)t_i.
$$
Since for each $i\leq j$ we have
$$(a+\epsilon_i)+(\epsilon_j+1-\epsilon_i)=a+\epsilon_j+1,$$
this inequality becomes
$$
   (a+\epsilon_j+1)\left(\binom{k+n-1}{n-1}-\sum_{i=0}^j t_i\right)<\sum_{i=j+1}^s (a+\epsilon_i)t_i.
$$
But $\sum_{i=0}^s t_i=\binom{k+n-1}{n-1}$, so the left-hand side equals $(a+\epsilon_j+1)\sum_{i=j+1}^s t_i$. Hence we obtain
$$(a+\epsilon_j+1)\sum_{i=j+1}^s t_i < \sum_{i=j+1}^s (a+\epsilon_i)t_i,$$
which is equivalent to
$$\sum_{i=j+1}^s (\epsilon_i-\epsilon_j-1)t_i>0.$$
This is exactly condition \eqref{eq:teza 5.1}, completing the proof.
\end{proof}
\begin{remark}
The numerical invariant appearing in Theorem~\ref{the:unexp} 
can be directly read off from the splitting type of the corresponding vector bundle, namely:
$$D_0^k(\A_Z)\big|_{\eL}\;\cong\;\O_{\eL}(-a)^{\oplus (N-\sum_{i=1}^s t_i)}\;\oplus\;\O_{\eL}(-(a+\epsilon_1))^{\oplus t_1}\;\oplus\;\O_{\eL}(-(a+\epsilon_2))^{\oplus t_2}\;\oplus\;\cdots\;\oplus\;\O_{\eL}(-(a+\epsilon_s))^{\oplus t_s}.$$
\end{remark}
\section{Examples}\label{sec:examples}

In this section we present two concrete examples together with links to external files and {\tt Singular} code~\cite{SingularCode}, which can be used to test further instances. We start with an example in $\P^4$ (followed by one in $\P^3$).

\begin{example}
Consider the set $Z$ of $25$ points in $\P^4$ listed in Table~\ref{tab:P4pts}. 
These coordinates arise as normal vectors to the hyperplanes of a crystallographic arrangement (see, e.g., \cite{Cuntz}).

\begin{table}[ht]
\centering
\small
\begin{tabular}{@{}llll@{}}
\toprule
\multicolumn{4}{c}{Point coordinates in $\P^4$} \\
\midrule
$(0:0:0:0:1)$ & $(0:0:0:1:0)$ & $(0:0:0:1:1)$ & $(0:0:1:0:0)$ \\
$(0:0:1:0:1)$ & $(0:0:1:1:0)$ & $(0:0:1:1:1)$ & $(0:1:0:0:0)$ \\
$(0:1:0:0:1)$ & $(0:1:0:1:1)$ & $(0:1:1:0:1)$ & $(0:1:1:1:1)$ \\
$(0:1:1:1:2)$ & $(1:0:0:0:0)$ & $(1:0:0:1:0)$ & $(1:0:0:1:1)$ \\
$(1:0:1:1:0)$ & $(1:0:1:1:1)$ & $(1:0:1:2:1)$ & $(1:1:0:1:1)$ \\
$(1:1:1:1:1)$ & $(1:1:1:1:2)$ & $(1:1:1:2:1)$ & $(1:1:1:2:2)$ \\
$(1:1:2:2:2)$ & & & \\
\bottomrule
\end{tabular}
\caption{The set $Z\subset \P^4$.}
\label{tab:P4pts}
\end{table}
A direct computation in {\tt Singular}~\cite{Singular} shows that the restricted bundle splits as
$$D_0(\A_Z)\big|_{\eL}\;\cong\;
   \O_{\eL}(-4)\;\oplus\;\O_{\eL}(-5)\;\oplus\;\O_{\eL}(-7)\;\oplus\;\O_{\eL}(-8),
$$
so the splitting type along $\eL$ is $(a_1,a_2,a_3,a_4)=(4,5,7,8)$. 
With the notation of Theorem~\ref{the:unexp} (here $k=1$), we may write this as
$$a=4,\qquad 
   \epsilon_0=0,\ \epsilon_1=1,\ \epsilon_2=3,\ \epsilon_3=4,\qquad 
   t_1=t_2=t_3=1.
$$
The criterion from Theorem~\ref{the:unexp} then yields:

\medskip
\noindent
\emph{Type $(a+\epsilon_0+k,\,a+\epsilon_0)=(5,4)$.}
$$0<\sum_{i=1}^3 t_i(\epsilon_i-\epsilon_0-1)= 1\cdot 0 +1\cdot 2+ 1\cdot 3 = 5,$$
so an unexpected hypersurface of type $(5,4)$ occurs.
Moreover, applying the general formula for the number of independent conditions imposed by a linear subspace $\H$ of dimension $n-2$ in $\P^n$ (see \cite[Lemma~2.1]{DHRST}), we obtain the following expression for the virtual dimension of the corresponding linear system:
\begin{equation}
\label{eq:vdimForm}
   \vdim([I(Z,\H)]_{a+\epsilon_j+k})=\binom{a+\epsilon_j+k+n}{n}-\sum_{i=0}^{a+\epsilon_j-1}(i+1)\binom{a+\epsilon_j+k-i+n-2}{n-2}-|Z|. 
\end{equation}
In the present case this gives
$$\binom{9}{4}-\sum_{i=0}^3 (i+1)\binom{7-i}{2}-25=126-105-25=-4,$$
while computer calculations show that the actual dimension equals 1. Hence we obtain the hypersurface $\mathcal{S}$, whose construction has been described in Section~4. This construction guarantees that the multiplicity of $\mathcal{S}$ in the points of $Z$ is at least~1, and the computations confirm that it is in fact exactly~1 at every point of $Z$. Note also that the ideal $I(Z)$ imposes independent conditions on quartic forms.

\smallskip
\noindent
\emph{Type $(a+\epsilon_1+k,\,a+\epsilon_1)=(6,5)$.}
$$0<\sum_{i=2}^3 t_i(\epsilon_i-\epsilon_1-1)= 1\cdot 1+1\cdot 2=3,$$
so an unexpected hypersurface of type $(6,5)$ also occurs. 
In this case the expected dimension is given by formula
$$\binom{10}{4}-\sum_{i=0}^4 (i+1)\binom{8-i}{2}-25=210-185-25=0,$$
Computer calculations show that the actual dimension equals $3$. As before, the construction ensures that the multiplicity in the points of $Z$ is at least~1, and computations confirm that it is exactly~1 at every point.

\smallskip
\noindent
\emph{Type $(a+\epsilon_2+k,\,a+\epsilon_2)=(8,7)$.}
$$0< t_3(\epsilon_3-\epsilon_2-1)= 1\cdot 0 =0,$$
so the condition fails. Direct computations in Singular shows that in this example 
$$\adim([I(Z,\H)]_8)=\vdim([I(Z,\H])_8=9,$$ 
thus no unexpected hypersurface of type $(8,7)$ arises.
\end{example}

\begin{remark}
All examples discussed in this paper, including the one above, arise from \emph{crystallographic arrangements} in the sense of Cuntz–Heckenberger. This class contains, in particular, arrangements associated with Weyl root systems; it is known that crystallographic arrangements are simplicial, though not every simplicial arrangement is crystallographic. During the preparation of this paper, we compiled a broader list of examples of this type, covering various dimensions and parameter values~$(n,k,|Z|,\adim,\vdim)$. For each configuration, the external file \cite{SingularCode} provides the defining coordinates, the exponents of the corresponding logarithmic bundle, and the actual and virtual dimensions of the associated systems.  This dataset illustrates several recurring patterns in the behaviour of unexpected hypersurfaces and can serve as a reference for further exploration.
\end{remark}

\begin{example}
Consider now the so–called Fermat (or Ceva) configuration of points in $\P^3$.
The set $Z$ consists of $31$ points whose homogeneous coordinates are listed in Table~\ref{tab:P3pts}. Here $e$ denotes a primitive third root of unity.

\begin{table}[ht]
\centering
\small
\begin{tabular}{@{}llll@{}}
\toprule
\multicolumn{4}{c}{Point coordinates in $\P^3$} \\
\midrule
$(1:0:0:0)$ & $(0:1:0:0)$ & $(0:0:1:0)$ & $(0:0:0:1)$ \\
$(1:1:1:1)$ & $(1:1:1:e)$ & $(1:1:1:e^2)$ & $(1:1:e:1)$ \\
$(1:1:e:e)$ & $(1:1:e:e^2)$ & $(1:1:e^2:1)$ & $(1:1:e^2:e)$ \\
$(1:1:e^2:e^2)$ & $(1:e:1:1)$ & $(1:e:1:e)$ & $(1:e:1:e^2)$ \\
$(1:e:e:1)$ & $(1:e:e:e)$ & $(1:e:e:e^2)$ & $(1:e:e^2:1)$ \\
$(1:e:e^2:e)$ & $(1:e:e^2:e^2)$ & $(1:e^2:1:1)$ & $(1:e^2:1:e)$ \\
$(1:e^2:1:e^2)$ & $(1:e^2:e:1)$ & $(1:e^2:e:e)$ & $(1:e^2:e:e^2)$ \\
$(1:e^2:e^2:1)$ & $(1:e^2:e^2:e)$ & $(1:e^2:e^2:e^2)$ & \\
\bottomrule
\end{tabular}
\caption{The Fermat (Ceva) configuration of points $Z\subset \P^3$.}
\label{tab:P3pts}
\end{table}

We now focus on the case $k=3$. Computations performed in \texttt{Singular} for the set $Z$ given in Table \ref{tab:P3pts} show that
$$D^3_0(\A_Z)\big|_{\eL}\;\cong\;\O_{\eL}^{\oplus 9}\;\oplus\;\O_{\eL}(-1)^{\oplus 3}\;\oplus\;\O_{\eL}(-2)^{\oplus 4}\;\oplus\;\O_{\eL}(-3)^{\oplus 2}\;\oplus\;\O_{\eL}(-4).$$
In the notation of Theorem \ref{the:unexp}, we therefore have
$$a=0,\qquad (\epsilon_i)_{i=1}^{4}=(1,2,3,4),\qquad (t_i)_{i=1}^{4}=(3,4,2,1).$$
Applying Theorem \ref{the:unexp} yields that there should exist an unexpected surface of type $(4,1)$, since the corresponding numerical inequality
$$ 0 < \sum_{i=2}^4 t_i(\epsilon_i-\epsilon_1-1)= 4\cdot 0+2\cdot 1 + 1\cdot 2=4$$
is satisfied.

However, in this case the construction from Theorem \ref{thm:syzConstr} gives rise to a surprising behaviour: one of the points of the configuration $Z$ becomes a \emph{triple point} of the resulting surface. More precisely, the point $(0:1:0:0)\in Z$ acquires multiplicity $3$. The subsequent computations must therefore take this phenomenon into account. Taking the triple point into account, we use the standard virtual–dimension count (as in \ref{eq:vdimForm}) augmented by the contribution of a triple point, which is $\binom{5}{3}$.  This gives
$$\vdim([I(Z,\H)]_4)= \binom{7}{3}-\binom{5}{1}-\binom{5}{3}-30=-10.$$
On the other hand, computations in \texttt{Singular} give $\adim([I(Z,\H)]_4)=3.$ Hence, the system exhibits a nontrivial excess dimension, indicating unexpected behaviour. Moreover, in this case the set $Z$ does not impose independent conditions on quartic forms. According to the \texttt{Singular} computation, four of the conditions coming from vanishing on $Z$ are redundant. Correcting the virtual dimension accordingly gives $\vdim([I(Z,\H)]_4)=-6,$ which still remains negative. Thus, the unexpectedness of the quartic surface arises from the occurrence of the triple point rather than from a simple failure of independence.

Finally, we remark that the appearance of multiple points within the fixed set $Z$, as a consequence of the syzygy-based construction, has already been pointed out in the literature (see \cite[Example 6.3]{MT-G}). In that planar case, two of the Fermat points became double, while here the three-dimensional Fermat configuration produces a \emph{triple} point. This suggests that the phenomenon persists in higher dimensions, reflecting the intrinsic combinatorial symmetry of the configuration.

\end{example}
\begin{remark}
    All computations discussed in the previous example were performed using a dedicated \texttt{Singular} script.  The code is freely available and can be used to reproduce the results presented here or to experiment with other configurations of points and multiplicities~\cite{SingularCode}. 
\end{remark}

	\subsection*{Acknowledgements.}

    Marek Janasz was partially supported by the National Science Centre (Poland) Sonata Bis Grant 2023/50/E/ST1/00025.

    Grzegorz Malara was partially supported by the National Science Centre (Poland) Sonata Grant 2023/51/D/ST1/00118.

	\footnotesize
	\noindent
	\textsc{Marek Janasz, Grzegorz Malara}: Department of Mathematics, University of the National Education Commission, 	Podchor\c a\.zych 2, 30-084 Krak\'ow, Poland, \\
    \textit{E-mail address:}  \texttt{marek.janasz@uken.krakow.pl}\\
    \textit{E-mail address:} \texttt{grzegorzmalara@gmail.com}\\
	\\

    \noindent
	\textsc{Halszka Tutaj-Gasi\'nska}: Faculty of Mathematics and Computer Science, Jagiellonian University, Stanis{\l}awa {\L}ojasiewicza 6, 30-348 Kraków, Poland,\\
	\textit{E-mail address:}  \texttt{halszka.tutajgasinska@gmail.com}
	
\end{document}